\documentclass{article}
\usepackage{amssymb}
\usepackage{amsmath}
\usepackage{amsthm}
\usepackage[usenames]{color}
\usepackage{tikz}

\usepackage{authblk} 

\newtheorem{statement}{}
\newtheorem{theorem}[statement]{Theorem}

\newtheorem{remark}[statement]{Remark}

\newcommand{\dd}[1]{{#1}^{\ast\ast}}

\usepackage{graphicx}
\usepackage{amssymb}
\usepackage{amsmath}
\usepackage[usenames]{color}
\usepackage{tikz}

\makeatletter
\newcommand{\subjclass}[2][1991]{%
  \let\@oldtitle\@title%
  \gdef\@title{\@oldtitle\footnotetext{#1 \textbf{Mathematics subject classification:} #2}}%
}
\newcommand{\keywords}[1]{%
  \let\@@oldtitle\@title%
  \gdef\@title{\@@oldtitle\footnotetext{\textbf{Keywords:} #1}}%
}
\makeatother

\def\N{\Bbb N}
\def\IN{\hbox{{\rm I}\kern-.13em{\rm N}}}

\def\IR{\hbox{{\rm I}\kern-.13em{\rm R}}}
\def\nin{n\in\N}

\def\ee{\varepsilon}

\def\dd{\delta}

\def\nn{\|\cdot\|}

\def\Qbar{\hbox{\sl Q\kern-.45em{\vrule height.63em width.05em
depth-.033em}}~}
\def\hrf{\hbox to.75in{\hrulefill}}

\newlength{\RoundedBoxWidth}
\newsavebox{\GrayRoundedBox}

\newenvironment{dbox}[1][\dimexpr\textwidth-2.1ex]%
   {\medskip\setlength{\RoundedBoxWidth}{\dimexpr#1}
    \begin{lrbox}{\GrayRoundedBox}
       \begin{minipage}{\RoundedBoxWidth}}%
   {   \end{minipage}
    \end{lrbox}
    \centerline{
    \begin{tikzpicture}%
       \draw node[draw=black,fill=black!5,rounded corners,%
             inner sep=1ex,text width=\RoundedBoxWidth]%
             {\usebox{\GrayRoundedBox}};
    \end{tikzpicture}
    }}

\newcommand{\caja}[1]{%
	\begin{dbox}{%
			{#1}}
\end{dbox}}

\begin{document}

\title{A remark on totally smooth renormings
}

\author{Ch. Cobollo\footnote{Christian Cobollo. 
Email: chcogo@upv.es 
}
,
A. J. Guirao\footnote{Antonio José Guirao. 
Email: anguisa2@mat.upv.es.
}
,
and V. Montesinos\footnote{Vicente Montesinos. 
Email: vmontesinos@mat.upv.es. Corresponding author.
}
}

\date{February 2020}

\affil{\textit{Instituto Universitario de Matem\'atica Pura y Aplicada. Universitat Polit\`ecnica de Val\`encia (Spain).}}




\maketitle

\begin{abstract}
E. Oja, T. Viil, and D. Werner showed, in {\em Totally smooth renormings}, Archiv der Mathematik, {\bf 112}, 3, (2019), 269--281, that a weakly compactly generated Banach space $(X,\nn)$ with the property that every linear functional on $X$ has a unique Hahn--Banach extension to the bidual $X^{**}$ (the so-called Phelps' property U in $X^{**}$, also known as the Hahn--Banach smoothness property) can be renormed to have the stronger property that for {\em every} subspace $Y$ of $X$, every linear functional on $Y$ has a unique Hahn--Banach extension to $X^{**}$ (the so-called total smoothness property of the space). We mention here that this result holds in full generality ---without any restriction on the space--- and in a stronger form, thanks to a result of M. Raja, {\em On dual locally uniformly rotund norms}, Israel Journal of Mathematics {\bf 129} (2002), 77--91.
\end{abstract}

\section{A totally smooth ---and more--- renorming}

The norm $\nn$ of a Banach space $(X,\nn)$ is said to be {\bf strictly convex} (or {\bf rotund}), if for $x,y\in S_X$ such that $\|x+y\|=2$ we have $x=y$. The norm is said to be {\bf locally uniformly convex} (or {\bf locally uniformly rotund}) ({\bf LUR}, for short), if $x\in S_X$, $x_n\in S_{X}$ for $\nin$, and $\|x+x_n\|\rightarrow 2$ implies $\|x_n-x\|\rightarrow 0$. The norm $\nn$ in a dual Banach space is said to have {\bf property $w^*$-LUR} if $x^*_n\rightarrow x^*_0$ in the $w^*$-topology as soon as $x^*_0,\ x^*_n\in S_{X^*}$ for $\nin$, and $\|x^*_n+x^*_0\|\rightarrow 2$.

\medskip

Attention has been paid to the problem of uniqueness of norm-preserving extensions (also called {\bf Hahn--Banach extensions}) of any continuous functional defined on a closed subspace $Y$ of a Banach space $X$ to the whole of $X$ (a property investigated by R. R. Phelps in \cite{phelps-u}). The norm $\nn$ of a Banach space $(X,\nn)$ has the so-called {\bf Hahn--Banach smooth property} ({\bf HBS}, for short) if every $x^*\in X^*$ has a unique norm-preserving  extension to $X^{**}$. A result that can be traced back to G. Godefroy \cite{godefroy}, and that appears also in \cite[Lemma III.2.14]{m-ideals}, says that this property is equivalent to the coincidence of the topologies $w$ and $w^*$ on the unit sphere $S_{X^*}$ of the dual space $X^*$ (we can call this property {\bf WW$^{*}$Kadets}, for short). A stronger property, called {\bf total smoothness} ({\bf TS} for short) is that for {\em any} closed subspace $Y$ of $X$, any $y^*\in Y^*$ has a unique Hahn--Banach extension to $X^{**}$. This is equivalent, by a result of A. E. Taylor and S. R. Foguel \cite{taylor}, \cite{foguel}, to the HBS property plus the rotundness of the dual norm $\nn^*$. In \cite{sullivan} it was proved that a {\em separable} space whose norm has the HBS property has a TS renorming, and in  \cite{oja} this result was extended to the class of {\em weakly compactly generated} spaces (i.e., spaces having a weakly compact linearly dense subset).

Here we just point out that the renorming result holds, even in a stronger form, {\em without any restriction on the space}. This observation is based on the following M. Raja's result: {\em If $X^*$ is a dual Banach space, then the two following conditions are equivalent: {\em (i)} $X^*$ admits an equivalent dual LUR norm, and {\em (ii)} $X^*$ admits an equivalent norm with the WW$^*$Kadets property} \cite{raja}.

\medskip

We believe that putting together those results as in Theorem \ref{thm-main} below may help to clarify the connections between the different properties mentioned above.

\medskip

\caja{
\begin{theorem}\label{thm-main}
Let $(X,\nn)$ be a Banach space. Then, the following statements are equivalent:

{\em (i)} $X$ has an equivalent norm with property HBS.

{\em (ii)} $X$ has an equivalent norm whose dual norm has property WW${^*}$Kadets.

{\em (iii)} $X$ has an equivalent norm whose dual norm is LUR.

{\em (iv)} $X$ has an equivalent norm with property TS.
\end{theorem}
}

\medskip
\begin{proof}

    (i)$\Leftrightarrow$(ii) is the aforementioned  result of Godefroy \cite{godefroy}. 
    
    (ii)$\Leftrightarrow$(iii) is the quoted result of Raja above \cite{raja}. Although not explicitly stated in this reference, it is simple to observe that the topology induced by a dual LUR norm $\nn^*$ coincides with the $w^*$-topology on the unit sphere $S_{X^*}$. Indeed, by the LUR property, given $x^*_0\in S_{X^*}$  and $\ee>0$, there exists $\dd>0$ such that if $x^*\in S_{X^*}$ satisfies $\|x^*_0+x^*\|>2(1-\dd)$, then $\|x^*_0-x^*\|<\ee$. Let $\{x^*_i\}$ be a net in $S_{X^*}$ that $w^*$-converges to $x^*_0$. Find, by Riesz's Lemma, $x_0\in S_{X}$ such that $\langle x_0,x^*_0\rangle>1-\dd$. There exists $i_0$ such that $\langle x_0,x^*_i\rangle>1-\dd$ for $i\ge i_0$. Thus, $\|x^*_0+x^*_i\|\ge \langle x_0,x^*_0+x^*_i\rangle>2(1-\dd)$ for $i\ge i_0$, hence $\|x^*_0-x^*_i\|<\ee$ for $i\ge i_0$, and the conclusion follows.

(iii)$\Rightarrow$(iv) follows from the Taylor--Foguel result \cite{taylor},\cite{foguel}, quoted above, and the observation in the proof of the equivalence (iii)$\Leftrightarrow$(ii) here.

(iv)$\Rightarrow$(i) is obvious.
\end{proof}

\begin{remark}\rm
\begin{enumerate}

\item Observe that, in particular, the TS norm defined in (iii) above on every Banach space with a HBS norm is Fr\'echet differentiable.

\item Banach spaces that satisfy one (and then all) of the conditions (i) to (iv) in Theorem \ref{thm-main} have been characterized in other different ways. Let us mention here that, for example, Theorem 1.4 in \cite{FOR} provides a few of them, in terms of (a) the existence of a dual norm in $X^*$ such that $(S_{X^*},w^*)$ is a Moore space, or (b) the existence of an equivalent dual norm such that $(S_{X^*},w^*)$ is symmetrized by a symmetric $\rho$ such that every point $x^*\in S_{X^*}$ has $w^*$-neighborhoods  of arbitrary small $\rho$-diameter, or (c) the existence of a dual equivalent norm such that $(S_{X^*},w^*)$ is metrizable, or even that (d) that $(B_{X^*},w^*)$ is a descriptive compact space (for details, see the op. cit. and the reference list there).
\end{enumerate}
\end{remark}

\begin{remark} \rm
Let us mention here (only with a hint for the proofs) that, for a Banach space $(X,\nn)$ whose norm $\nn$ has property HBS,

\begin{enumerate}
\item The norm, restricted to any closed subspace of $X$, has property HBS too, a consequence of the $w^*$-lower semicontinuity of the dual norm.

\item $X$ is Asplund, as it follows from (i) above and a separable reduction argument.

\item  $X$ is {\bf nicely smooth} (i.e., there is no proper $1$-norming subspace in $X^*$) and that, in fact, every James boundary is strong (see, e.g., \cite[Paragraph 3.11.8.3]{y2}).

\item If $(X,\nn)=(C(K),\nn_{\infty})$, where $K$ is a compact topological space, then $K$ is finite. This follows from the fact that the set of extreme points of $B_{C(K)^*}$ is $\{\pm\dd_{k}:\ k\in K\}$, that all extreme points are distributed between two closed hyperplanes, the Krein--Milman theorem, and the consequent reflexivity of the space $C(K)$. This observation depends strongly on  the fact that the norm on $C(K)$ is the supremum norm. A space $C(K)$, for $K$ an infinite compact space, may admit an equivalent norm $\nn$ whose dual is LUR (and so $\nn$ has property HBS): Just take $K$ an infinite countable compact space; it is metrizable and scattered (see, e.g., \cite[Lemma 14.21]{y2}), hence $C(K)$ is Asplund (see, e.g., \cite[Theorem 14.25]{y2}). Thus, $C(K)^*$ is separable, and the conclusion follows from a classical result of Kadets (see, e.g., \cite[Section 2]{fmz}).

\item There exists a LUR renorming of $X$. This follows from the aforementioned Raja's result and a result of R. Haydon in  \cite{haydon}. Note that it is an open problem (see, e.g.,  \cite[Problem 1]{smith-troyanski} and \cite[Problem 102]{open}) whether a space $X$ has a LUR renorming as soon as it has a norm whose dual norm has property $w^*$-LUR. \end{enumerate}
\end{remark}

\section*{Acknowledgements}
Supported by  AEI/FEDER (project MTM2017-83262-C2-2-P of Ministerio de Econom\'{\i}a y Competitividad), by Fundaci\'on S\'eneca, Regi\'on de Murcia (grant 19368/PI/14), and Universitat Polit\`ecnica de Val\`encia (A. J. Guirao)

Supported by AEI/FEDER (project MTM2017-83262-C2-1-P of Ministerio de Econom\'{\i}a y Competitividad) and Universitat Polit\`ecnica de Val\`encia (V. Montesinos)

\section*{Acknowledgements}


\begin{thebibliography}{}


\bibitem{y2} M. Fabian, P. Habala, P. H\'ajek, V. Montesinos, and V. Zizler, {\em Banach Space Theory: The Basis of Linear and Nonlinear Analysis}, Springer, New York, Dordrecht, Heidelberg, London (2011)
    
\bibitem{fmz} M. Fabian, V. Montesinos, and V. Zizler, {\em Smoothness in Banach spaces. Selected problems}, Rev. R. Acad. Cien. Serie A. Mat. RACSAM, {\bf 100} (1-2), 101--125 (2006) 

\bibitem{FOR} S. Ferrari, J. Orihuela, M. Raja, {\em Generalized metric properties of spheres and renorming of Banach spaces}. Rev. R. Acad. Cienc. Exactas Fís. Natl. Ser. A Math. RACSAM {\bf 113}, 2655--2663 (2019)

\bibitem{foguel} S. R. Foguel, {\em On a theorem by A. E. Taylor}, Proc. Amer. Math. Soc. {\bf 9}, 325 (1958)

\bibitem{godefroy} G. Godefroy, {\em Points de Namioka, espaces normants, applications \`a la th\'eorie isom\'etrique de la dualit\'e}, Israel J. Math. {\bf 38}, 209--220 (1981)

    \bibitem{open} A. J. Guirao, V. Montesinos, and V. Zizler, {\em Open Problems in the Geometry and Analysis of Banach Spaces}, Springer International Pub. Switzerland (2016)

\bibitem{m-ideals} P. Harmand, D. Werner, and W. Werner, {\em M-ideals in Banach Spaces and Banach Algebras}, Lecture Notes in Math. {\bf 1547}, Springer, Berlin (1993)

\bibitem{haydon}  R. Haydon, {\em Locally uniformly rotund norms in Banach spaces and their duals}. J. Funct. Anal. {\bf 254},  2023--2039 (2008)

\bibitem{oja} E. Oja, T. Viil, and D. Werner, {\em Totally smooth renormings}, Archiv der Mathematik, {\bf 112}, 3,  269--281 (2019)

\bibitem{phelps-u} R. R. Phelps, {\em Uniqueness of Hahn–Banach extensions and unique best approximation}, Trans. Amer. Math. Soc. {\bf 95}, 238--255 (1960)

\bibitem{raja} M. Raja, {\em On dual locally uniformly rotund norms}, Israel Journal of Mathematics {\bf 129}, 77--91 (2002)

\bibitem{smith-troyanski} R. J. Smith, S. L. Troyanski, {\em Renormings of $C(K)$ spaces}. Rev. R. Acad. Cienc. Exactas Fís. Natl. Ser. A Math. RACSAM {\bf 104} (2), 375-–412 (2010)

\bibitem{sullivan} F. Sullivan, {\em Geometrical properties determined by the higher duals of a Banach space}, Illinois J. Math. {\bf 21}, 315--331 (1977)

\bibitem{taylor} A. E. Taylor, {\em The extension of linear functionals}, Duke Math. J. {\bf 5}, 538--547 (1939)


\end{thebibliography}
\end{document}